\documentclass[12pt]{iopart}

\usepackage{iopams}  

\usepackage{epsfig,psfrag,epsf}
\begin{document}

\title[A regularized weighted least gradient problem]{A regularized weighted least gradient problem for conductivity imaging}
\author{Alexandru Tamasan$^1$
 and Alex Timonov$^2$}
\address{$^1$ Department of Mathematics, University of Central Florida, Orlando, FL, USA}
\address{$^2$ Division of Mathematics and Computer Science, University of South Carolina Upstate, Spartanburg, SC, USA}
\ead{tamasan@math.usf.edu, atimonov@uscupstate.edu}
\vspace{10pt}
\begin{indented}
\item[]March 2018
\end{indented}

\begin{abstract}
We propose and study a regularization method for recovering an approximate electrical conductivity solely from the magnitude of one interior current density field. Without some minimal knowledge of the boundary voltage potential, the problem has been recently shown to have nonunique solutions, thus recovering the exact conductivity is impossible. The method is based on solving a weighted least gradient problem in the subspace of functions of bounded variations with square integrable traces. The computational effectiveness of this method is demonstrated in numerical experiments.
\end{abstract}

\section{Introduction}
\label{intro}

Let $\Omega\subset R^n$, $n=2,3$, be a Lipschitz domain modeling a conductive body. We revisit the inverse hybrid problem of reconstructing an inhomogeneous, isotropic, electrical conductivity $\sigma$ from knowledge of the magnitude of one current density field inside $\Omega$. The problem may be reduced to solving a singular, degenerate elliptic equation (the 1-Laplacian in a conformal Euclidean metric) subject to various boundary conditions \cite{seo,NTT07}, or can be cast as  a minimization problem involving a weighted gradient term \cite{NTT09,MNT,NTV16}. Without some minimal knowledge of the voltage potential at the boundary, the problem has non-unique solution as recently characterized in \cite{NTV16}; where additional measurements of the voltage potential along a curve joining the electrodes were proposed to establish uniqueness. Other approaches, some of which are mentioned below, assume  knowledge of the magnitude of two current density fields, or of the entire field. The only known modality of obtaining the interior data involves rotations in a magnetic resonance machine \cite{joy}. This makes any boundary voltage potential measurement, while not impossible, at least impractical. 

In this paper we address the non-uniqueness via a regularization method, which recovers an approximate conductivity without recourse to any boundary voltage information. In any vicinity of the given interior data, we identify some ``ideal"  data, which uniquely determines the sought conductivity. By reversing the roles, and interpreting the available data as a perturbation of the ideal one, we then propose a reconstruction method, and analyze the continuous dependence problem. Numerical experiments will show feasibility of the method.

The forward problem is modeled by the Robin boundary conditions, which, in the case of two electrodes, we show to be equivalent to the Complete Electrode Model \cite{somersaloCheneyIsaacson}. More precisely, assume that a current density field is generated by injecting/extracting  a net current $I>0$ from a couple of surface electrodes $e_{\pm}$ assumed  bounded Lipschitz subdomains in $\partial\Omega$, with real valued impedance $z>0$. For a known conductivity $\sigma$, the voltage potential $u_0\in H^1(\Omega)$ distributes inside according to
\begin{eqnarray}
& & \nabla\cdot\sigma\nabla u_0=0, \quad\mbox{in}\;\Omega,\label{conductivity_eq}\\
& & \sigma\frac{\partial u_0}{\partial\nu}=-b_0 u_0+c_0,\quad \mbox{on}\; \partial\Omega,\label{robin} 
\end{eqnarray}
where
\begin{equation}\label{coefficients}
b_0:=\frac{1}{z}\left\{
\begin{array}{l}
1 \mbox{ on }e_\pm,\\
0, \mbox{ off } e_{\pm}, 
\end{array}\right.
\qquad
\mbox{and}\qquad
c_0:=\left\{
\begin{array}{l}
\pm I, \mbox{ on }e_\pm,\\
0, \mbox{ off }e_{\pm},
\end{array}\right.
\end{equation}and $\nu$ denotes the outer unit normal to the boundary, 

By replacing the conductivity in \ref{conductivity_eq} by $a/|\nabla u_0|$, the problem reduces to solving a boundary value problem for a generalized 1-Laplacian as originally proposed in \cite{seo}. The work in \cite{NTT07} was first to point out the connection with minimum surfaces in a Riemannian space determined by the interior data, and proposed a method to recover the conductivity from Cauchy data. For Dirichlet data in \cite{NTT09, NTT10} the problem was reduced to  minimum gradient problem for functions of given trace at the boundary, and, in \cite{NTV16}, extended to the Complete Electrode Model (CEM) boundary conditions originally introduced in \cite{somersaloCheneyIsaacson} . Existence and/or uniqueness of such weighted gradient problems were studied in \cite{jerrard99} and \cite{MNT17}, with extensions to perfectly insulated and conducting inclusions in \cite{MNTa_SIAM,MNT17}. A structural stability result for the minimization problem can be found in \cite{nashedTa10}. Reconstruction algorithms based on the minimization problem were proposed in \cite{NTT09} and \cite{MNT}, and based on level set methods in \cite{NTT07,NTT09,verasTamasanTimonov14}. Continuous dependence on $\sigma$ on $a$ (for a given unperturbed Dirichlet data) can be found in \cite{montaltoStefanov}, and, for partial data in \cite{montaltoTamasan}.  For further references on determining the isotropic conductivity based on measurements of current densities see \cite{zhang,seo,seo_ieee, jeun-rock, nachman,lee,  kimLee}, and for reconstructions on anisotropic conductivities from multiple measurements  see \cite{MaDeMonteNachmanElsaidJoy,MNNick2014, balGuoMonard2d, balGuoMonard3d}.

In here we seek to determine an approximate conductivity $\sigma$, solely from knowledge of the magnitude
\begin{equation}\label{a0}
a_0:=|\sigma \nabla u_0|
\end{equation}
of the current density field inside $\Omega$,
where $u_0\in H^1(\Omega)$ (functions and their gradient are square integrable) is the unique solution to the Robin problem \ref{conductivity_eq} and \ref{robin}.

As recently characterized in \cite{NTV16}, we note that $\sigma$ is not uniquely determined by $a_0$. For example, for
any $\varphi: Range(u_0)\to Range(u_0)$ an increasing Lipschitz continuous function, satisfying
$\varphi(t)=t$ for $t\in u_0(e_+)\cup u_0(e_-)$, one can verify that $u_\varphi=\varphi\circ u_0$ is another solution of the Robin problem corresponding to the conductivity $\sigma/(\varphi'\circ u_0)$, while the magnitude of the induces current density field does not change.

Using the original idea in \cite{NTT09}, we approach the inverse problem via a weighted minimum gradient problem, here modeled for Robin boundary conditions. 

For some nonegative $a\in C(\Omega)\cap L^\infty(\Omega)$  (playing the role of an $L^2$-approximation of the given data $a_0$),   $b\in L^\infty(\partial\Omega)$, and $h$ harmonic function to be specified later,  we  consider the minimization of the functional
\begin{equation}\label{WeightFuncEpsBV}
v\mapsto G(v;a) :=\int_\Omega a|\nabla v|dx + \frac{1}{2}\int_{\partial\Omega} b (v-h)^2 ds.
\end{equation}
With one exception it will suffice to minimize the functional  over $H^1(\Omega)$. However, for the continuous dependance result in Theorem \ref{main_thm} we need to consider the functional over the subspace $BV_2(\Omega)$ of functions of bounded variation with square integrable traces. This is the smallest subspace, in which a minimizing sequence is compact. In this regard, when $v\in BV_2(\Omega)$, the first integral term will be understood in the sense of a Radon measure $|D v|$  applied to a bounded continuous $a$. More precisely,
\begin{equation}\label{BVa}
| Dv|(a):=\sup \{\int_{\Omega}v\nabla \cdot F dx:\: F\in C^1_0(\Omega; R^d)),\;|F(x)|\leq a(x)\}.
\end{equation}

In Section 2 we will provide a triple of coefficients $(a,b,h)$, such that the functional  $G(\cdot;a)$ will satisfy the following existence and uniqueness property:
\begin{equation}\label{U}
\mbox{There exists } u\in H^1(\Omega) \mbox{ the unique minimizer of } G(\cdot;a) \mbox{ over }BV_2(\Omega).
\end{equation}
Our interior data $a_0$ in \ref{a0} may not be bounded in the vicinity of the boundary of the electrodes. Moreover,  while $u_0$ will minimize the functional $G(\cdot; a_0)$ in $H^1(\Omega)$, it will not be unique, since, for any $\varphi$ as in the counterexample above,  $\varphi\circ u_0$ will also be a minimizer. This motivates us to consider  the regularized functional
\begin{equation}\label{RegWeightFuncEpsBV}
G^\delta(v;a):=\int_\Omega a|\nabla v|dx+ \frac{1}{2}\int_{\partial\Omega} b (v-h)^2 ds+\frac{\delta}{2}\int_\Omega|\nabla (v-h)|^2dx.
\end{equation}
The continuous dependence of the minimizing sequence with respect to the weight $a$ in $L^2(\Omega)$, and $\delta\to 0$ (studied in Section 4) constitute the basis of the numerical method used in Section 5. 

To connect with the work in \cite{NTV16},  we remark here that,  for two electrodes, the Robin problem is equivalent to the Complete Electrode Model (CEM) problem in  \cite{somersaloCheneyIsaacson}, up to a scaling factor.  Not essential, but simplifying the exposition, we further assume the electrodes have equal surface areas, 
$\displaystyle |e|:= |e_\pm|$.
In the complete electrode model 
the voltage potential $v$ solution of \ref{conductivity_eq} inside $\Omega$, and an unknown constant voltage $V$ satisfy the boundary conditions
\begin{eqnarray}
&v+z\sigma\frac{\partial v}{\partial\nu}= \pm V,\quad \mbox{on}\; e_\pm,\label{robin_kForward} \\
&\int_{e_\pm} \sigma\frac{\partial v}{\partial\nu}ds = \pm I \label{inject_kForward},\\
&\frac{\partial v}{\partial\nu}=0, \quad\mbox{on}\; \partial\Omega\setminus(e_+\cup e_-). \label{no_flux_off_electrodesForward}
\end{eqnarray}
Under the assumptions that  $\Omega$ is a Lipschitz domain, $\sigma$ is essentially bounded away from zero and infinity, the electrodes $e_\pm$ have positive impedance and are (relatively) open connected subsets of $\partial\Omega$ with disjoint closure, the CEM problem has a unique solution
$(v; V)\in H^1(\Omega)\times R$, see \cite{somersaloCheneyIsaacson}, or the appendix in \cite{NTV16}.  For an arbitrary $\lambda>0$ the pair $(\lambda u_0,\lambda z I)$ clearly solves 
\ref{conductivity_eq}, \ref{robin_kForward} (with $V=z\lambda I$), and \ref{no_flux_off_electrodesForward}.
An application of Green's theorem in the Robin model yields $\displaystyle\int_{e_+}u_0ds=-\int_{e_-}u_0ds,$ which well defines the scaling choice
\begin{equation*}
\lambda^{-1}:=\left(|e| -\frac{1}{zI}\int_{e_+}u_0ds\right)=\left( |e| +\frac{1}{zI}\int_{e_-}u_0ds\right).
\end{equation*}
With this choice of scaling, one can check that $\lambda u_0$ also satisfies \ref{inject_kForward}, and thus  $\lambda u_0=v,\mbox{ and } a_0=\lambda \sigma|\nabla v|.$ Therefore, if we use to magnitude of the current density field corresponding to the Robin problem or to the CEM,  we would recover the same conductivity $\displaystyle \sigma= a_0/|\nabla u_0|=|\sigma\nabla v|/|\nabla v|.$

\section{Remarks on the smoothness of solutions to the Robin problem}
\label{smoth}

Our techniques, which is based on the minimization of the functional \ref{WeightFuncEpsBV}, requires the weight $a$ be bounded continuous in $\Omega$. This regularity cannot be achieved solely on the smoothness in the conductivity $\sigma$, as the regularity of the coefficients appearing in the Robin condition \ref{robin} also play a role. Throughout we assume a conductivity
 \begin{equation}\label{sigma}
 \sigma\in C^{1/2}(\overline\Omega)\mbox{ with } \sigma|_{\partial\Omega}\in C^2(\partial\Omega).
 \end{equation}
Under this smoothness assumption, the elliptic regularity for solutions to the Robin problem  (e.g., \cite[Theorem 7.4, Remark 7.2]{LionsMagenes}) yields that $u_0\in C^{1/2}(\overline\Omega)\cap C^{1,1/2}(\Omega)$. Moreover, $\nabla u_0$, and thus, $a_0$  extend by H\"{o}lder-continuity to all points in $\partial\Omega\setminus\partial e_\pm$, see  \cite[Proposition B.1. (ii)]{NTV16} for details. However, the right hand side of \ref{robin} is merely in $H^{1/2-s}$ for some $s>0$, which yields $u_0 \in H^{3-s}(\Omega)$, insufficient to conclude the boundedness of $\nabla u_0$ in three dimensions. Namely, at the boundary of the electrodes,  the tangential derivative normal to  $\partial e_\pm$  may blow up, yielding an unbounded interior data $a_0$ in \ref{a0}. However,  if we considered some $C^2$- smooth approximations of $b_0$ and $c_0$ that made the right hand side of \ref{robin} lie in $H^{1/2}(\partial\Omega)$, then the same bootstrap argument in the proof of \cite[Proposition B.1. (ii)]{NTV16} would apply to show that the corresponding Robin solution $u \in H^3(\Omega)$. Indeed,  for a right hand side of \ref{robin} in $H^{1/2}(\partial\Omega)$, the solution $u\in H^2(\Omega)$, which in turn yields $u\in H^{3/2}(\partial\Omega)$, which together with $C^2$-smoothness of the coefficients, yield that the right hand side of \ref{robin} now lie in  $H^{3/2}(\partial\Omega)$. Another application of the classical regularity result yields $u\in H^3(\Omega)\subset C^{1,1/2}(\overline\Omega)$. Thus, in two and three dimensions, $a=\sigma|\nabla u| \in C^{1/2}(\overline\Omega)$ is bounded continuous.

\section{Existence and uniqueness of a minimizer}
\label{ex_un}

The regularized method can be better understood through a family of forward problems.  For each $\epsilon>0$ small, let first define the boundary function 
\begin{equation*}
\tilde{b}_\epsilon:=\frac{1}{z}\left\{
\begin{array}{l}
1 \mbox{ on }e_\pm,\\
\epsilon, \mbox{ on } \partial\Omega\setminus (e_-\cup e_+), 
\end{array}\right.
\end{equation*}
and recall the coefficients $b_0$ and $c_0$ in (\ref{coefficients}). According to the regularity remark in Section 2, the solution to the problem (\ref{conductivity_eq}) subject to a boundary condition as in (\ref{robin}) with $\tilde{b}_\epsilon$ replacing $b_0$, might not be of bounded gradient as needed. 

This  motivates to further consider, for each $\epsilon>0$, some smoother approximates (e.g., by gluing) $b_\epsilon\in C^2 (\partial\Omega)$ of $\tilde{b}_\epsilon$, respectively $c_\epsilon\in C^2(\partial\Omega)$ of $c_0$, with the only necessary property that 
\begin{equation} \label{lim_e}
\lim_{\epsilon\to 0}\|{b}_\epsilon-b_0\|_\infty=0, \mbox{ and }
\lim_{\epsilon\to 0}\|{c}_\epsilon-c_0\|_\infty=0.
\end{equation}

Let  $u_\epsilon\in H^1(\Omega)$ be the solution of the Robin problem  (\ref{conductivity_eq}) subject to
\begin{equation}\label{robinBC}
\sigma\frac{\partial u}{\partial\nu}=-b_\epsilon u+c_\epsilon, \quad\mbox{ on }\partial\Omega,
\end{equation}
and define an ``ideal" interior data $a_\epsilon$ as the magnitude of the corresponding current density field 
\begin{equation}\label{aepsilon}
a_\epsilon:=|\sigma \nabla u_\epsilon|.
\end{equation}
The remark in the section above shows that $a_\epsilon\in C^{1/2}(\overline\Omega)$ for $\epsilon>0$. Moreover, classical arguments on the continuous dependence (in particular, since the coercivity constant is bounded below independently of $\epsilon$), also apply to yield
\[
\|u_\epsilon-u_0\|_{H^1(\Omega)}\to 0,\mbox{ and } \|a_\epsilon-a_0\|_{L^2(\Omega)}\to 0, \mbox{ as }\epsilon\to 0^+.
\]
For each $\epsilon\geq 0$ small, it is convenient to consider  the harmonic function $h_\epsilon$, solution to 
\begin{equation}\label{harmonic}
\Delta h_\epsilon=0,\mbox{ in }\Omega,\quad h_\epsilon|_{\partial\Omega}=\frac{c_\epsilon}{b_\epsilon},
\end{equation}
where $b_\epsilon$ and $c_\epsilon$ 
are as introduced above.

For each $\epsilon\geq 0$ small, let us consider  the functional in (\ref{WeightFuncEpsBV}) corresponding to $a_\epsilon$, $b_\epsilon$, and $h_\epsilon$
\begin{equation}\label{1LaplFuncEps}
G (v;a_\epsilon):=\int_\Omega a_\epsilon |\nabla v|dx + \frac{1}{2}\int_{\partial\Omega} b_\epsilon (v-h_\epsilon)^2 ds,
\end{equation} 
and recall that for $\epsilon>0$ the functional extends over functions in
\begin{equation}\label{BV2}
BV_2(\Omega):=\{u:BV(\Omega):~ u|_{\partial\Omega}\in L^2(\partial\Omega)\},
\end{equation}
where the first integral is in the sense of the Radon measure $|Dv|$ applied to $a_\epsilon$ as in (\ref{BVa}).

The following result shows the regularizing effect of $\epsilon>0$.

\noindent
{\bf Theorem 1.}
{\it Let $\sigma$ satisfy (\ref{sigma}),  $b_0$, $c_0$, $u_0$,  and $a_0$ be as above and $G(\cdot;a_0)$ be as in (\ref{1LaplFuncEps}) with $\epsilon=0$.Then 
\begin{equation}\label{globalMIN}
G(u_0;a_0)\leq G(v;a_0), \quad\mbox{for all } v\in H^1(\Omega).
\end{equation}
Moreover, for $\epsilon>0$,  let $b_\epsilon$, $c_\epsilon$, $u_\epsilon$,  $a_\epsilon$ and $G(\cdot;a_\epsilon)$ be as above in (\ref{1LaplFuncEps}). Then $u_\epsilon\in C^{1,1/2}(\overline\Omega)$ is the unique minimizer of $G(\cdot;a_\epsilon)$  in $BV_2(\Omega)$,
\begin{equation}\label{minProblem}
u_\epsilon=\mbox{argmin}\{ G(v;a_\epsilon):\; v\in BV_2(\Omega)\}.
\end{equation}
In particular, the exact conductivity can be recovered uniquely from $a_\epsilon$ by  
\begin{equation}\label{sigmaEps}
\sigma=\frac{a_\epsilon}{|\nabla u_\epsilon |}.
\end{equation}
}

{\bf Proof.} For any $v\in H^1(\Omega)$, we estimate
\begin{eqnarray}
G(v;a_\epsilon)&= \int_\Omega a_\epsilon |\nabla v|dx +\frac{1}{2}\int_{\partial\Omega}b_\epsilon (v-h_\epsilon)^2 ds\nonumber\\
&=\int_\Omega \sigma |\nabla u_\epsilon| |\nabla v|dx +\frac{1}{2}\int_{\partial\Omega}b_\epsilon (v-h_\epsilon)^2 ds\nonumber\\
&\geq\int_\Omega\sigma \nabla u_\epsilon\cdot\nabla vdx+\frac{1}{2}\int_{\partial\Omega}b_\epsilon (v-h_\epsilon)^2 ds\nonumber\\
&=\int_{\partial\Omega}(-b_\epsilon u_\epsilon+c)vds +\frac{1}{2}\int_{\partial\Omega}b_\epsilon (v-h)^2 ds\nonumber\\
&=\frac{1}{2}\int_{\partial\Omega}b_\epsilon(v-u_\epsilon)^2ds+
\frac{1}{2}\int_{\partial\Omega}b_\epsilon(h_\epsilon^2-u_\epsilon^2)ds\nonumber\\
&\geq \frac{1}{2}\int_{\partial\Omega}b_\epsilon(h_\epsilon^2-u_\epsilon^2)ds =G(u_\epsilon;a_\epsilon),\label{estimates}
\end{eqnarray}
where the second equality uses (\ref{aepsilon}), the third equality uses the divergence theorem and the fact the $u_\epsilon$ solves the Robin problem (\ref{conductivity_eq}), (\ref{robinBC}). This proves (\ref{globalMIN}). 

We show next that $u_\epsilon$ is a global minimizer of the functional over the larger set $BV_2(\Omega)$. Let $v\in BV_2(\Omega)$ be arbitrary. 
By mollification (e.g.,  see \cite[Remark 2.12]{giusti}), there exists a sequence $\displaystyle\{v_n\}\subset W^{1,1}(\Omega)$ with ${v_n}|_{\partial\Omega}=v|_{\partial\Omega}$, and such that $v_n\to v$ in $L^1(\Omega)$, and
\begin{equation}\label{convInMeas}
\lim_{n\to\infty}\int_\Omega a_\epsilon|\nabla v_n|dx = |Dv|(a_\epsilon). 
\end{equation}
By taking the limit with $n\to\infty$  in 
\[
G_\epsilon(u_\epsilon;a_\epsilon)\leq \int_\Omega a_\epsilon |\nabla v_n|dx+\int_{\partial\Omega}b_\epsilon(v_n-h_\epsilon)^2ds=\int_\Omega a_\epsilon |\nabla v_n|dx+\int_{\partial\Omega}b_\epsilon(v-h)^2ds,
\]
and using  (\ref{convInMeas}), we conclude that
$\displaystyle G_\epsilon(u_\epsilon;a_\epsilon)\leq|D v|(a_\epsilon)+\int_{\partial\Omega}b_\epsilon(v-h_\epsilon)^2ds.$

Now let $v\in BV_2(\Omega)$ be another minimizer of $G (\cdot; a_\epsilon)$, and  consider a mollified sequence $\{v_n\}\subset W^{1,1}(\Omega)$,  $v_n|_{\partial\Omega}=v|_{\partial\Omega}$ as above (\cite[Remark 2.12]{giusti}) to estimate
\begin{eqnarray*}
G(v;a_\epsilon)&= |Dv |( a_\epsilon)+\frac{1}{2}\int_{\partial\Omega}b_\epsilon (v-h)^2 ds\nonumber\\
&=\lim_{n\to\infty}\int_\Omega \sigma |\nabla u_\epsilon| |\nabla v_n|dx +\frac{1}{2}\int_{\partial\Omega}b_\epsilon (v-h)^2 ds\nonumber\\
&\geq \limsup_{n\to\infty}\int_\Omega\sigma \nabla u_\epsilon\cdot\nabla v_ndx+\frac{1}{2}\int_{\partial\Omega}b_\epsilon (v-h)^2 ds\nonumber\\
&=\limsup_{n\to\infty}\int_{\partial\Omega}(-b_\epsilon u_\epsilon+c)v_nds +\frac{1}{2}\int_{\partial\Omega}b_\epsilon (v-h)^2 ds\nonumber\\
&=\int_{\partial\Omega}(-b_\epsilon u_\epsilon+c)vds +\frac{1}{2}\int_{\partial\Omega}b_\epsilon (v-h)^2 ds\nonumber\\
&=\frac{1}{2}\int_{\partial\Omega}b_\epsilon(v-u_\epsilon)^2ds+
\frac{1}{2}\int_{\partial\Omega}b_\epsilon(h^2-u_\epsilon^2)ds\nonumber\\
&\geq \frac{1}{2}\int_{\partial\Omega}b_\epsilon(h^2-u_\epsilon^2)ds =G_\epsilon(u_\epsilon;a_\epsilon),
\end{eqnarray*}
Since $v$ is also a minimizer $G_\epsilon(v;a_\epsilon)=G_\epsilon(u_\epsilon;a_\epsilon)$ and all the inequalities above hold with equality, in particular
\[
\int_{\partial\Omega}b_\epsilon(v-u_\epsilon)^2ds=0.
\]

For $\epsilon>0$, the weight $b_\epsilon$ is essentially positive on the boundary,  which yields
\begin{equation}\label{Dirichlet}
v |_{\partial\Omega}= u_\epsilon|_{\partial\Omega}.
\end{equation}
Next, we note  that, for competitors  restricted to the affine subspace 
$$
D_\epsilon:=\{v\in BV_2(\Omega):\;v|_{\partial\Omega}=u_\epsilon|_{\partial\Omega}\},
$$ the minimization problem
$\min\left\{G_\epsilon (v;a_\epsilon):\; v\in D_\epsilon \right\}$
is equivalent to
\begin{equation}\label{DirichletProb}
\min\left\{ |Dv|( a_\epsilon): \quad v\in BV(\Omega), \; v|_{\partial\Omega}=u_\epsilon|_{\partial\Omega} \right\}.
\end{equation}
Since $u_\epsilon\in C^{1,1/2}(\overline\Omega)$ is a solution, we apply the uniqueness result \cite[Theorem 1.1]{MNT17} to the minimization problem (\ref{DirichletProb}) to conclude that
$$v=u_\epsilon, \quad \mbox{ in } \Omega.$$

Following from the definition of $a_\epsilon$ in (\ref{aepsilon}), and the strict positivity of the conductivity, the set of critical points $\{x\in \Omega:\;|\nabla u_\epsilon|=0\}$ coincide with the set of zeros of $a_\epsilon$. Since the set of critical points is negligible  in $\Omega$, the equality (\ref{sigmaEps}) holds almost everywhere. Since $\sigma$ is assumed continuous, the equality (\ref{sigmaEps}) must then hold at all points in $\Omega$.

\section{Regularization of the weighted least gradient problem}
\label{regl}

Since our available data is not $a_\epsilon$ but rather  the $L^2$ approximate $a_0$, we cannot apply Theorem 1  directly to recover $\sigma$. Moreover the functional 
$\displaystyle v\mapsto \int_{\Omega}a_0|\nabla v|dx + \int_{\partial\Omega}b_\epsilon(v-h_\epsilon)^2ds
$ may not have a minimizer (the previous arguments based on the forward problem no longer work, since we mix the internal data $a_0$ coming from $b_0,c_0$ in (\ref{robin}), with the regularized coefficients $b_\epsilon$ and $c_\epsilon$). This motivates us to consider the regularized functional below, where, for brevity, we drop the $\epsilon$-subscript from the notations.

For some nonnegative $a\in L^2(\Omega)$, $b\in L^\infty(\partial\Omega)$ with
\[
0<\frac{\epsilon}{z}\leq b\leq\frac{1}{z}, \mbox{ a.e. } \partial\Omega,
\]
and $h$ harmonic in $\Omega$, consider the
functional
\begin{equation}\label{RegWeightFunct}
G^\delta(v;a):=\int_\Omega a|\nabla v|dx+ \frac{1}{2}\int_{\partial\Omega} b (v-h)^2 ds+\frac{\delta}{2}\int_\Omega|\nabla (v-h)|^2dx.
\end{equation}
The sum of the quadratic terms 
\begin{equation}\label{goodF}
F^\delta(v):=\frac{1}{2}\int_{\partial\Omega}b v^2ds +\frac{\delta}{2}\int_{\Omega}|\nabla v|^2 dx
\end{equation}  
in (\ref{RegWeightFunct}) gives an equivalent (square of the) norm in $H^1(\Omega)$, since
\begin{equation}
\min\left\{\frac{\epsilon}{2z},\frac{\delta}{2}\right\} \|u\|_1^2
\leq F^\delta(u)\leq \max\left\{\frac{1}{2z},\frac{\delta}{2}\right
\} \|u\|_1^2,\label{equivNorm}
\end{equation}where
$$\|u\|^2_1:=\int_{\partial\Omega}|u|^2ds+\int_{\Omega}|\nabla u|^2dx.$$

The unique minimizer in $H^1(\Omega)$ of the functional (\ref{RegWeightFunct}) follows from classical convex minimization arguments, which we include them below for completeness.

{\bf Proposition 1.}~{\it 
For $\epsilon,\delta>0$ arbitrarily fixed, and ${a}\in L^2(\Omega)$ positive,  let $G^\delta(\cdot,{a})$ be the functional in (\ref{RegWeightFunct}). The minimization problem
\[
\min\{ G^\delta (v; {a}): \; v\in H^1(\Omega) \}
\]
has a unique solution.
}

{\bf Proof.} We show first that  $G^\delta(\cdot;{a})$ is weakly lower semi-continuous. Let 
$v_n\rightharpoonup v$ be a weakly convergent sequence in $ H^1(\Omega)$. we need to show that
\begin{equation}\label{lsc}
G^\delta(v;{a})\leq \liminf_{n\to\infty} G^\delta(v_n;{a}).
\end{equation}
The weak lower semicontinuity of $F^\delta_\epsilon:H^1(\Omega)\to R$
follows from a classical argument that uses its convexity
\[
F^\delta_\epsilon(v_n)\geq F^\delta_\epsilon (v) + \int_{\partial\Omega}b_\epsilon v(v_n-v)ds+\delta\int_\Omega\nabla v\cdot\nabla(v_n-v)dx,
\]
and Fatou's lemma.

The weak lower semicontinuity of the weighted gradient functional
\begin{equation}\label{weightGrad}
v\mapsto \int_{\partial \Omega}{a}|\nabla v|dx
\end{equation}
uses some classical arguments in the theory of functions of bounded variation:
Let $\{a_m\}$ be an increasing sequence of bounded continuous
functions, which converges in $L^2(\Omega)$ sense to ${a}$. For each fixed index $m$, let $f=(f_1,...,f_n)\in
C_0^1(\Omega;R^n)$ be arbitrary with $|f|\leq a_m$. Since
$v_n\rightharpoonup v$ in $L^2(\Omega)$ we have
\begin{eqnarray}\label{estimate_1}
&\int_\Omega v \nabla\cdot fdx =\lim_{n\to\infty}\int_\Omega v_n \nabla\cdot fdx
=\liminf_{n\to\infty}\int_\Omega v_n\nabla\cdot fdx\nonumber\\
&\leq\liminf_{n\to\infty}~\sup\left\{\int_\Omega v_n\nabla\cdot
gdx:~~g\in C^1_0(\Omega; R^n),~~|g|\leq a_m\right\}\nonumber\\
&=\liminf_{n\to\infty}\int_{\Omega}a_m|\nabla v_n|dx\leq
\liminf_{n\to\infty}\int_{\Omega}{a}|\nabla v_n|dx,
\end{eqnarray}
where the last inequality above uses the fact that $a_m\leq {a}$. By taking the supremum in (\ref{estimate_1}) over all $f\in
C_0^1(\Omega;R^n)$ with $|f|\leq a_m$ we get
\begin{eqnarray}\label{estimate_2}
\int_{\Omega}a_m|\nabla
v|dx&=\sup\left\{\int_\Omega v\nabla\cdot fdx:~~f\in
C^1_0(\Omega; R^n),~~|f|\leq a_m\right\}\nonumber\\
&\leq\liminf_{n\to\infty}\int_{\Omega}a_0|\nabla v_n|dx.
\end{eqnarray}
By letting $m\to\infty$ in (\ref{estimate_2}) we obtain the weakly lower 
semi-continuity for (\ref{weightGrad}).

We showed that $v\mapsto G^\delta(v;{a}):H^1(\Omega)\to R$ is weakly lower semi-continuous in $H^1(\Omega)$. Since 
$G^\delta(\cdot;{a})$ is also strictly convex, it has a unique minimizer. 

\section{Convergence properties of the regularized minimizing sequence}
\label{mainresultsection}

For $a$ bounded continuous in $\Omega$ satisfying
\begin{equation}\label{positivity}
\inf(a)=:\alpha>0,
\end{equation}
$b\in L^\infty(\partial\Omega)$ with
\begin{equation}\label{lower_b}
0<\frac{\epsilon}{z}\leq b\leq\frac{1}{z}, \mbox{ a.e. } \partial\Omega,
\end{equation}
and $h\in H^1(\Omega)$ harmonic in $\Omega$, recall the functional (\ref{WeightFuncEpsBV})
\begin{equation}\label{Recall:WeightFuncEpsBV}
G(v;a)=|Dv(a)|+\int_{\partial\Omega}b(v-h)^2ds.
\end{equation}

In this section we assume that $G(\cdot;a)$ satisfies the existence and uniqueness hypothesis (\ref{U}), and propose a minimization scheme based on the regularized problem in Section \ref{regl}.

Note that Theorem 1 yields that (\ref{U}) holds for $a=a_\epsilon$ (and $b=b_\epsilon$, and $h=h_\epsilon$).  

We will often use the trivial identity 
\begin{equation}\label{useful}
\int_\Omega \tilde{a}|\nabla v|dx= \int_{\Omega}a|\nabla v|dx+\int_{\Omega}(\tilde{a}-a) |\nabla v|dx,
\end{equation}
which allows us to exchange two arbitrary weights $ a,\tilde{a}\in L^2(\Omega)$. For brevity we use $\|\cdot\|$ to denote the $L^2(\Omega)$-norm, and by $\|\cdot\|_1$ the $H^1(\Omega)$-norm.

{\bf Theorem 2.}~{\it
Let $\Omega\subset R^d$ be a bounded Lipschitz domain with connected boundary, and $a\in BC(\Omega)$ satisfy (\ref{positivity}). Assume that the functional $G(\cdot;a)$ in (\ref{Recall:WeightFuncEpsBV}) satisfy the hypothesis (\ref{U}), and let $u\in H^1(\Omega)$ be the unique minimizer
\begin{equation}\label{minimum}
u=\mbox{\it argmin}\{G(v;a):\;v \in BV_2(\Omega)\}.
\end{equation}
 Let $\{a_n\}\subset L^2(\Omega)$ be a sequence of positive functions, with 
\begin{equation}\label{L2_depend}
\|a_n-a\|\longrightarrow 0, \mbox{ as } n\to\infty.
\end{equation}
and $\delta_n\downarrow 0$ be a decreasing
sequence such that
\begin{equation}\label{not_so_fast}
\lim_{n\to\infty}\frac{\|a_n-a\|^2}{\delta_n}=0.
\end{equation}
Corresponding to each $n$, consider the regularized functional  $v\mapsto G^{\delta_n}(v, a_n)$ as in (\ref{RegWeightFunct}), and let
\begin{equation}\label{min_sequence}
u_n:= \mbox{\it argmin} \{ G^{\delta_n}(v, a_n):\;v\in H^1(\Omega)\}
\end{equation}
be the corresponding minimizer provided by Proposition 1. Then
\begin{equation}
\lim_{n\to\infty}G^{\delta_n}_\epsilon (u_n;a_n) =\lim_{n\to\infty}G(u_n;a)=G(u;a).\label{minimizing_sequence}
\end{equation}
Moreover, on a subsequence $\{\tilde{u}_n\}$ of $\{u_n\}$,
\[
\tilde{u}_n\longrightarrow u, \mbox{ in } L^q(\Omega),\;0\leq q\leq\frac{d}{d-1},
\]
and, for any open subset $O\subset\Omega$,
\begin{equation}\label{onpatches}
\lim_{n\to\infty}\int_{O}a_n|\nabla \tilde{u}_n|dx=\lim_{n\to\infty}\int_{O}a|\nabla \tilde{u}_n|dx =\int_{O}a |\nabla u|dx. 
\end{equation}
}

{\bf Proof.}
Despite the fact that $\|u_n\|_1$ may not be uniformly bounded, we prove first that
\begin{equation}\label{educatedecay}
\lim_{n\to\infty}\int_\Omega (a-a_n)|\nabla u_n|dx=0.
\end{equation}
Let  $n$ be sufficiently large so that $\|a_n\|\leq 2\|{a}\|$. Recall the functional $F^{\delta_n}(\cdot; a_n)$ in (\ref{goodF}) with $\delta=\delta_n$ and $a=a_n$, and the induced norm on $H^1(\Omega)$ in (\ref{equivNorm}). We estimate
\begin{eqnarray}
\min\left\{\frac{\epsilon}{2z},\frac{\delta_n}{2}\right\} \|u_n-h\|_1^2&\leq F^{\delta_n} (u_n - h;a_n)\nonumber\\
&\leq  F^{\delta_n} (u_n - h;a_n)+\int_{\Omega}a_n |\nabla u_n|dx\nonumber\\
&=G^{\delta_n} (u_n; a_n)  \nonumber\\
&\leq G^{\delta_n}(h; a_n) \nonumber \\
&=\int_\Omega a_n|\nabla h|dx\nonumber\\
&\leq 2\|{a}\| \|\nabla h\|,\label{boundEst}
\end{eqnarray}
where the third inequality uses the minimizing property defining $u_n$. Note that the right hand side of (\ref{boundEst}) is  independent of $\delta_n$ to yield:
\begin{equation}\label{uniformBoundDelta}
\|u_n\|_1 \leq {C}\max\left\{\frac{2z}{\epsilon},\frac{2}{\delta_n}\right\}^{1/2},
\end{equation}
for some constant $C$ dependent on $\|{a}\|$ and $\|h\|_1$. In particular since $\delta_n\to 0$, for sufficiently large $n$, we obtained,
\begin{equation}\label{u_nbounded}
\|u_n\|_1\leq C\frac{1}{\sqrt{\delta_n}}, 
\end{equation}
where $C$ depends only on $\|{a}\|$ and the $\|h\|_1$. The rate of decay (\ref{not_so_fast}) together with (\ref{u_nbounded}) yields (\ref{educatedecay}).

Since ${u}$ is a minimizer of $G(\cdot ;{a})$, we estimate
\begin{eqnarray}
G ({u}; {a})&\leq \liminf_{n\to\infty}G (u_n; {a})\leq \limsup_{n\to\infty}G (u_n; {a})\nonumber\\
&=\limsup_{n\to\infty}\left\{G (u_n;a_n)+\int_\Omega({a}-a_n)|\nabla u_n|dx\right\}\nonumber\\
&\leq\limsup_{n\to\infty}\left\{G (u_n;a_n)+\frac{\delta_n}{2}\int_\Omega|\nabla( u_n-h)|^2dx +\int_\Omega({a}-a_n)|\nabla u_n|dx \right\}\nonumber\\
&=\limsup_{n\to\infty}\left\{G^{\delta_n}(u_n;a_n)+\int_\Omega({a}-a_n)|\nabla u_n|dx\right\}\nonumber\\
&=\limsup_{n\to\infty}G^{\delta_n}(u_n;a_n).\label{onedirection}
\end{eqnarray}
where the first equality uses (\ref{useful}), the next to the last equality uses the definition of $G^{\delta_n}_\epsilon$, and the last equality uses (\ref{educatedecay}). Similarly,
\begin{eqnarray}
G ({u}; {a})&\leq \liminf_{n\to\infty}G (u_n; {a})\nonumber\\
&=\liminf_{n\to\infty}\left\{G (u_n;a_n)+\int_\Omega({a}-a_n)|\nabla u_n|dx\right\}\nonumber\\
&\leq\liminf_{n\to\infty}\left\{G (u_n;a_n)+\frac{\delta_n}{2}\int_\Omega|\nabla( u_n-h)|^2dx +\int_\Omega({a}-a_n)|\nabla u_n|dx \right\}\nonumber\\
&=\liminf_{n\to\infty}\left\{G^{\delta_n} (u_n;a_n)+\int_\Omega({a}-a_n)|\nabla u_n|dx\right\}\nonumber\\
&=\liminf_{n\to\infty}G^{\delta_n}(u_n;a_n)\leq\limsup_{n\to\infty}G^{\delta_n} (u_n;a_n).
\label{onedirection2}
\end{eqnarray}

The reverse inequality also holds
\begin{eqnarray}
&\limsup_{n\to\infty}G^{\delta_n}( u_n ;a_n)\leq\limsup_{n\to\infty}G^{\delta_n}( {u} ;a_n)\nonumber\\
&=\limsup_{n\to\infty}\left\{G^{\delta_n}({u} ;{a})+\int_\Omega(a_n-{a})|\nabla {u}|dx\right\}\nonumber\\
&=\limsup_{n\to\infty}\left\{G(u ;{a})+\frac{\delta_n}{2}\int_\Omega|\nabla( u-h)|^2dx+\int_\Omega(a_n-{a})|\nabla u|dx\right\}\nonumber\\
&=G(u;{a}).\label{otherdirection}
\end{eqnarray}
In the estimate (\ref{otherdirection}), the first inequality uses (\ref{min_sequence}),  while the last equality uses (\ref{L2_depend}) and the assumption $u\in H^1(\Omega)$. 

The inequalities (\ref{onedirection}), (\ref{onedirection2}), and (\ref{otherdirection}) prove the identity (\ref{minimizing_sequence}). In particular, we showed that
\begin{equation}
\lim_{n\to\infty}\left(\int_\Omega a|\nabla u_n|dx+\int_{\partial\Omega}b(u_n-h)^2ds\right)=\int_\Omega a|\nabla u|dx+\int_{\partial\Omega}b(u-h)^2ds.\label{minimizing_sequence_2}
\end{equation}

Note that both the regularization
parameter and the coefficients in the functional are changing with $n$. In particular, the sequence $u_n$ may not be
bounded in $H^1(\Omega)$.  However, we show next that $\{u_n\}$ is bounded in
$W^{1,1}(\Omega)$; endowed with the norm
$$\|u\|_{1,1}:=\int_{\partial\Omega}u ds+\int_\Omega|\nabla u|dx.$$
Recall the lower bound $\alpha$ in (\ref{positivity})  to estimate
\begin{eqnarray}
&\min\{\alpha,2\sqrt{b}\}\|u_n\|_{1,1}\leq\int_{\Omega}{a}|\nabla u_n|dx+
\int_{\partial\Omega}2\sqrt{b}|u_n|ds\nonumber\\
&\leq\int_{\Omega}{a}|\nabla u_n|dx+\int_{\partial\Omega}2\sqrt{b_\epsilon}|u_n-h|ds+\int_{\partial\Omega}2\sqrt{b}|h|ds\nonumber
\\
&\leq\int_{\Omega}{a}|\nabla u_n|dx+\int_{\partial\Omega}{b_\epsilon}(u_n-h)^2ds+ |\partial\Omega|+
\int_{\partial\Omega}2\sqrt{b}|h|ds\nonumber\\
&\leq G(u_n;{a})+ |\partial\Omega|+\int_{\partial\Omega}2\sqrt{b}|h|ds\nonumber\\
&=G(u_n;a_n)+\int_\Omega ({a}-a_n)|\nabla u_n|dx+ |\partial\Omega|+\int_{\partial\Omega}2\sqrt{b}|h|ds\nonumber\\
\leq &G(h ; a_n)+ C\frac{\|{a}-a_n\|} {\sqrt{\delta_n}}+ |\partial\Omega|+\int_{\partial\Omega}2\sqrt{b}|h|ds\nonumber\\
&=G(h ; {a})+ \int_\Omega (a_n-{a})|\nabla h|dx+ C\frac{\|{a}-a_n\|} {\sqrt{\delta_n}}+ |\partial\Omega|+\int_{\partial\Omega}2\sqrt{b}|h|ds\nonumber\\
&=G(h ; {a})+ \|a_n-{a}\|\|\nabla h\|+ C\frac{\|{a}-a_n\|} {\sqrt{\delta_n}}+ |\partial\Omega|+\int_{\partial\Omega}2\sqrt{b}|h|ds,\label{1,1bound}
\end{eqnarray}
where the fifth inequality uses the bound (\ref{u_nbounded}). By the hypothesis (\ref{not_so_fast}) on the rate of decay of $\delta_n$, the right hand side above is uniformly bounded in $n$. 

Since  $\min\{\alpha,\mbox{\it essinf}(b)\}>0$ we  showed that $\|u_n\|_{1,1}$ is uniformly bounded. An application of Rellich-Kondrachov's compactness embedding
(e.g., \cite{ziemer}) shows the existence of a convergent
subsequence $\{\tilde{u}_n|_{\partial\Omega}\}$, with $\tilde{u}_n\to u^*$ in $L^q(\Omega)$ for all
$1\leq q< d/(d-1)$.  Moreover, since $\{\tilde{u}_n\}$ is bounded in $W^{1,1}(\Omega)$, the limit $u^*\in BV(\Omega)$ and $u^*|_{\partial\Omega}\in L^1(\partial\Omega)$, see e.g.,\cite{giusti}. 

Also following  from the estimate (\ref{1,1bound}), the sequence of traces  $\{u_n|_{\partial\Omega}\}$ is uniformly bounded in $L^2(\partial\Omega)$. In particular, $u^*|_{\partial\Omega}\in L^2(\partial\Omega)$, and, possibly passing to a subsubsequence, $\{\tilde{u}_n\}$ converges weakly in $L^2(\partial\Omega)$ to $u^*|_{\partial\Omega}$. 

We show next that  $\{\tilde{u}_n\}$ also converges strongly in $L^2(\partial\Omega)$. We recall the weak lower semi-continuity properties on each of the two functionals in $G$. The first one is the lower semi-continuity of the total variations. For any ${a} \in C(\Omega)\cap L^\infty(\Omega)$,  
\begin{equation}\label{lscDu}
|Du^*| (a)\leq\liminf_{n\to\infty} \int_\Omega{a} |\nabla \tilde{u}_n | dx.
\end{equation}
The second is the weak lower semi-continuity of the quadratic term,
\begin{equation}\label{lscsqr}
\frac{1}{2}\int_{\partial\Omega}b(u^*-h)^2ds\leq\liminf_{n\to\infty}\frac{1}{2}\int_{\partial\Omega}b(u_n-h)^2ds.
\end{equation}
By adding (\ref{lscDu}) and (\ref{lscsqr}) and using (\ref{minimizing_sequence_2}) we get
\begin{equation}\label{mustbeequality} 
G(u^*,a)\leq\liminf_{n\to\infty} G(\tilde{u}_n;a)=G(u;a).
\end{equation}
Since $u$ was assumed the unique minimizer of $G$ in $BV_2(\Omega)$, we conclude
that equality must hold in (\ref{mustbeequality}), i.e., $\displaystyle G(u^*;a)=\liminf_{n\to\infty}G(\tilde{u}_n;a)$, and that
\begin{equation}\label{uniq}
u=u^*. 
\end{equation}
Moreover, each of the inequalities  (\ref{lscDu}) and (\ref{lscsqr}) must also be equalities. By possibly passing to a further sub-subsequence, we have shown that
\begin{eqnarray}
&|Du|(a)=\lim_{n\to\infty}\int_{\Omega}a|\nabla\tilde{u}_n|dx, \nonumber \\
&\int_{\partial\Omega}b(u-h)^2ds=\lim_{n\to\infty}\int_{\partial\Omega} b(\tilde{u}_n-h)^2ds.\label{convL2norm}
\end{eqnarray}
For any $O$ an open subset of $\Omega$, the arguments of \cite[Theorem 5.2.3]{ziemer} carries verbatim to conclude the upper semi-continuity property, for a
\[
\limsup_{n\to\infty}\int_{O}a|\nabla \tilde{u}_n|dx \leq\int_O a|\nabla u|dx.
\]
Now (\ref{onpatches}) follows by an application of (\ref{educatedecay}).

Note that the convergence in $L^2(\partial \Omega)$-norm in (\ref{convL2norm}) together with the weak convergence yield the strong convergence $\tilde{u}_n|_{\partial\Omega}\to u|_{\partial\Omega}$ in $L^2(\partial\Omega)$.

\section{Numerical demonstration}
\label{num}

We demonstrate the computational feasibility of the regularized method in some numerical experiments. The detailed numerical study will be presented elsewhere.

\begin{figure}[ht]
\centerline{\vbox{\hbox{\includegraphics[width=16cm]{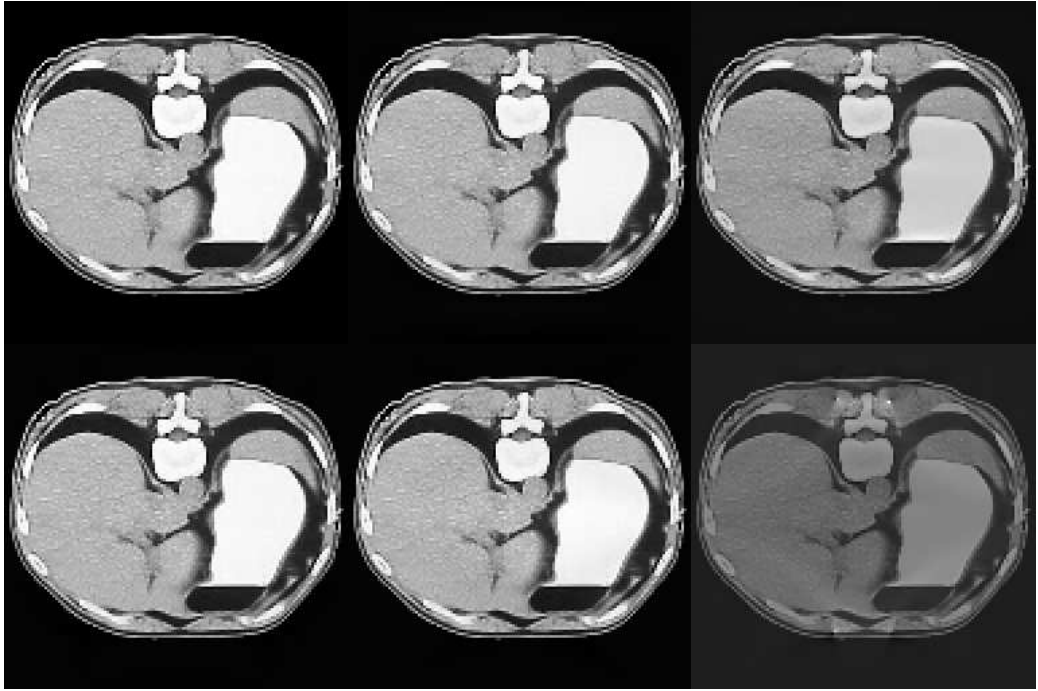}}}} 
\caption{Comparison of the reconstructed mean conductivity distributions. The level of the roundoff and truncation errors in the interior data does not exceed $10^{-5}$. The parameters $\epsilon = 5\cdot 10^{-4}$ and $\delta = 3\cdot 10^{-3}$ are chosen.}
\label{figure_1}\end{figure}

In all the numerical experiments below, the original conductivity $\sigma$ is simulated on a real abdominal CT image (shown in the left upper corner in Figure \ref{figure_1}) of a human. The image is embedded into a unit square $S=[0,1]\times[0,1]$, so that the space between the image and sides of the square is filled with a homogeneous medium with $\sigma(x) = 1$. The image is rescaled to the realistic range [1, 1.8] S/m of the electrical conductivity typical to the biological tissues. The interior data, i.e., the magnitude of the current density $a = \sigma|\nabla{u}|$, is simulated by solving the forward problem for the conductivity equation with the CEM. We use the standard Galerkin finite element method for computing its numerical solution. 

To solve the regularized minimization problem (\ref{RegWeightFunct}) in Section \ref{regl}, we use an iterative procedure based on solving forward Robin problems for updated conductivities. While similar to the algorithm developed in \cite{NTT09} in connection to the Dirichlet problem, at each iteration, we now solve 
\begin{eqnarray}
& & \nabla\cdot\left(\sigma + \delta\right)\nabla{u} = 0~\mbox{in}~\Omega, \label{cond_reg} \\
& & \left(\sigma + \delta\right)\frac{\partial{u}}{\partial{\nu}} + b_{\epsilon}u = \delta\frac{\partial{h}}{\partial{\nu}}~\mbox{on}~\partial\Omega. \label{rob_reg}
\end{eqnarray}

We simulate the experiments using electrodes of two apertures. In the full aperture case, the electrodes span the entire top respectively bottom side. In the smaller aperture case, the electrodes are centered and span one half of the upper/lower side. The reconstructed images obtained by the new method are also compared to those obtained by the alternating split Bregman algorithm proposed and developed in \cite{MNT} for the Dirichlet problem. For the latter, we use the calculated trace of the Robin solution as the needed Dirichlet data in the alternating split Bregman algorithm.

Numerically, we use the finite differences approach, where the resulting linear system is solved by an implicit conjugate gradient method, in which the preconditioned matrix is inverted and the correction vector is computed on a Krylov subspace in each iteration. All computations were performed on the Dell Precision workstation T5400 running under IDL 6.2.

%
%
Figure~\ref{figure_1} demonstrates comparison of the original conductivity distribution (shown in the left upper corner) with the conductivity means recovered from the interior data. In the upper row we show the conductivity means obtained from the interior data simulated for the full electrode apertures, i.e., the electrode length coincides with the size length. The conductivity mean obtained by the proposed algorithm is shown in the middle of the upper row. Its relative error is $4\cdot 10^{-3}$, whereas the relative error of the conductivity mean obtained by the alternating split Bregman algorithm shown in the right upper corner is $1.5\cdot 10^{-2}$. In the lower row we show the conductivity means for the reduced electrode apertures: half aperture (the left corner), two step sizes $2h$ (the middle - the proposed algorithm, the right corner - the alternating split Bregman algorithm). The corresponding relative $l_{2}$-errors of reconstruction are $3\cdot 10^{-3}$, $6\cdot 10^{-3}$, and $3\cdot 10^{-2}$, respectively.

\section{Conclusions}
We recover an approximate electrical conductivity from knowledge of the magnitude of the current density inside, without any knowledge of a voltage potential at the boundary. The new method relies on a solving a minimum weighted gradient problem corresponding to some Robin boundary conditions, which is regularized to mitigate for the elliptic degeneracy present in the problem. A compactness property of the minimizing sequence is shown in the space of functions of bounded total variation. Numerical experiments are conducted to demonstrate the feasibility of the method, and they are also compared to one of the method that uses full knowledge of the voltage potential at the boundary.

\pagebreak

\noindent
{\bf References}
\\

\end{document}